\documentclass{article}
\usepackage{amsmath}
\usepackage[affil-it]{authblk}

\usepackage[margin=2.0cm]{geometry}

\usepackage[normalem]{ulem}
\usepackage{amsmath,amssymb,mathtools,epsfig,epstopdf,mathpazo,textcomp,graphicx,chemarrow,soul}
\usepackage[dvipsnames]{xcolor}
\usepackage[toc,page]{appendix}
\usepackage[colorlinks]{hyperref}
\hypersetup{
    linkcolor=red,     
    citecolor=blue,    
    filecolor=cyan,    
    urlcolor=magenta   
}

\usepackage[nameinlink,capitalise]{cleveref}
\crefformat{equation}{(#2#1#3)}

\usepackage{cite}

\newcommand{\be}{\begin{eqnarray}}
\newcommand{\ee}{\end{eqnarray}}

\usepackage{xr}
\graphicspath{{./figSet/labelFigures/}}

\begin{document}

\title{Novel Generic Models for Differentiating Stem Cells\\ Reveal Oscillatory Mechanisms}

\def\correspondingauthor{\footnote{Corresponding author: a.rocco@surrey.ac.uk}}

\author[1]{Saeed Farjami}
\author[2]{Karen Camargo Sosa}
\author[3]{\\ Jonathan H.P. Dawes}
\author[2]{Robert N. Kelsh}
\author[1,4]{Andrea Rocco \correspondingauthor}

\affil[1]{Department of Microbial Sciences, University of Surrey, GU2 7XH, Guildford, UK }
\affil[2]{Department of Biology and Biochemistry, University of Bath, BA2 7AY, Bath, UK }
\affil[3]{Department of Mathematical Sciences, University of Bath, BA2 7AY, Bath, UK }
\affil[4]{Department of Physics, University of Surrey, GU2 7XH, Guildford, UK }

\date{}

\maketitle

\begin{abstract}
\noindent Understanding cell fate selection remains a central challenge in developmental biology. We present a class of simple yet biologically-motivated mathematical models for cell differentiation that generically generate oscillations and hence suggest alternatives to the standard framework based on Waddington's epigenetic landscape. The models allow us to suggest two generic dynamical scenarios that describe the differentiation process. In the first scenario gradual variation of a single control parameter is responsible for both entering and exiting the oscillatory regime. In the second scenario two control parameters vary: one responsible for entering, and the other for exiting the oscillatory regime. We analyse the standard repressilator and four variants of it and show the dynamical behaviours associated with each scenario. We present a thorough analysis of the associated bifurcations and argue that gene regulatory networks with these repressilator-like characteristics are promising candidates to describe cell fate selection through an oscillatory process.
\end{abstract}


\section{Introduction}
\label{sec:intro}

In 1940, Waddington proposed representing the complex regulatory dynamics driving the process of cellular differentiation as an 'epigenetic landscape' \cite{Waddington40} through which a single cell can be thought of as travelling. The differentiating cell encounters successive `decision points' in the landscape morphology that correspond to differentiation events. These decision points emerge from dynamical changes in the underlying cellular gene regulatory network (GRN), and are mathematically described by variations in parameters that influence the GRN dynamics. For example in zebrafish it is known that Wnt signalling plays a fundamental role in specification and commitment of melanocytes \cite{Vibert16,Dorsky00,Lewis04}, and the level of Wnt expression may therefore be considered as a parameter, driving the cellular GRN through one or more of the decision points hypothesised by Waddington.

Despite the philosophical attractiveness of the epigenetic landscape metaphor, the details have remained unclear and no completely self-consistent mathematical description has emerged. For example there are debates regarding the types of bifurcations characterising differentiation: saddle-node bifurcations \cite{Ferrel12,Huang17} have been proposed as being more consistent with the expected topologies of the core GRN responsible for differentiation than pitchfork bifurcations. Experimental support for these mathematical bifurcation phenomena, e.g. in embryonic stem cells (ESCs), is also lacking. For example, when progesterone is washed away from mature Xenopus oocytes, they do not dedifferentiate but remain mature \cite{Jin01}, a finding that is more consistent with the saddle-node bifurcation than with Waddington's landscape.

Of course, the detailed biological understanding of cell differentiation dynamics has grown significantly recently and it is appreciated that qualitatively different mechanisms may be at work in different contexts. For instance, detailed study of the differentiation of ESCs into two cell types that form either embryonic or extraembryonic tissues indicates that this does not occur in a simple deterministic manner \cite{Huang10}. Instead, cells appear to pass transiently through intermediate states, each seemingly primed to differentiate into one of the possible cell fates \cite{Wigger17,Graf08}. Hence ESCs might maintain a collection of meta-stable states, and a better framework to understand and interpret experimental data is required \cite{Graf08}.

Most fundamentally, the GRN topologies generating such dynamic stem cell states are poorly understood. In the present work we focus on the question of how broader multipotency (i.e. more than two fates) might generically be generated, in a way that allows differentiation to be biased towards locally-favourable outcomes. We note this is a topic of much current interest in systems biology, particularly noting the exploration of connections between GRN topology and oscillatory dynamics \cite{Perez18}, generic links with differentiation dynamics \cite{Suzuki11,Rabajante15,Kaity18}, and e.g. reshaping the epigenetic landscape by rewiring the GRN using cross-repression among transcription factors (TFs) with self-regulation \cite{Rabajante15}.

In this paper, we analyse and compare five minimally-constructed GRNs that admit cyclical dynamics and allow evolution towards fate commitment controlled by an external stimulus. Our GRNs are variants of the standard repressilator, a simple but fundamental circuit made up of three genes repressing each other. Theoretical modelling of this circuit and its subsequent engineering in {\em E. coli} have shown its capability to reproduce oscillatory behaviours under broad parameter ranges \cite{Elowitz00}.

We are here interested in adopting this model and the variants proposed to mimick entry into the dynamical oscillatory state characterizing multipotency, and the exit from it, towards a fully differentiated cell state.  

Therefore, besides exhibiting oscillations, the system must be equipped with a mechanism to exit the oscillatory regime and proceed towards differentiation. Hence the system must exhibit a sequence of at least two bifurcations, accounting for entering and exiting the oscillatory phase.

We identify two generic scenarios for entering and exiting the oscillatory regime, and we test both in our proposed GRN models. In the first scenario (S1), we hypothesise that a single parameter (i.e. one intervening signalling pathway) drives the system both into and out of the oscillatory regime. An example is Wnt signalling, required for the induction of neural crest, but also for specification of both melanocyte and sensory neuron fates through activation of key TFs (Mitf and Neurogenin, respectively) \cite{Dorsky00,Lewis04,Saint-Jeannet97,Dorsky98,Dunn00, Takeda00,Jin01,Lee04}, which we envisage as having this two-fold role. First, increasing Wnt signalling promotes neural crest cells to enter an oscillatory multipotent phase in which fate specification TFs are cyclically expressed; secondly, at higher signalling levels, it promotes the cell's exit from the oscillatory phase allowing specification of single cell fates through increased Mitf or Neurogenin expression which drive melanocyte or sensory neuron differentiation respectively.

In the second scenario (S2) we hypothesise again that a first bifurcation to oscillatory behaviour  is driven by an external signalling pathway, but its influence is then blocked above a critical value. Exit from the oscillatory regime results from other signalling pathways, associated with other parameters. This would arise for a multipotent progenitor that in response to a specific signal starts oscillating, as in S1, driving expression of multiple fate specification TFs, each characteristic of a different 'sub-state'.  Close to each transcriptional sub-state, the stem cell would express distinct cell signalling receptors; then sufficient activation of these receptors by environmental ligands could provide a mechanism that would move the cell out of the oscillatory phase, thus driving differentiation.

The paper is structured as follows. In \cref{sec:repressilator}, we consider the simplest GRN exhibiting oscillatory behaviour, the so-called repressilator \cite{Elowitz00}. We find that the system exhibits a parameter regime where trajectories tend towards an attracting limit cycle featuring a slowdown when passing near three sub-states; this appears as a rather weak effect but is numerically detectable. However, exit from the oscillatory regime leads to a single equilibrium, making this GRN unable to select between alternative differentiated states.

This result motivates \cref{sec:model1,sec:model2} where we extend the standard repressilator to include a second repressive circuit opposing the first; we therefore term this GRN the `cross-repressilator'. Within this new GRN we consider whether the two TFs act on each gene as either an `OR gate' or an `AND gate', leading to two variants of this circuit. We find that for the `OR gate', the circuit first transitions to oscillatory dynamics, as the standard repressilator, but then followed by a further transition to a stably differentiated cell state, in which only one of the genes is stably expressed. The `AND gate' is even more interesting, since it allows simultaneous expression of two out of the three genes. In both cases our results appear to hold for a large class of modelling choices for these gene interactions.

In \cref{sec:four_genes} we extend the `AND gate' cross-repressilator circuit to four and five genes. We show that co-expression of up to three genes is possible depending on the topology of the GRN. These dynamics are compatible with both S1 and S2, as they provide an exit mechanism dependent on a single or multiple bifurcation parameters into a number of alternative equilibria. Furthermore, we note that in all the models, when a cell is in the oscillatory phase, scoring for gene expression of fate-specific TFs would result (in a snapshot view) in the cell being considered `fate specified', despite the fact that it actually retains full multipotency.

A discussion  of the biological implications of our models is presented in \cref{sec:conc}.

\section{The standard repressilator}
\label{sec:repressilator}

\begin{figure}[ht]
\begin{center}
\includegraphics{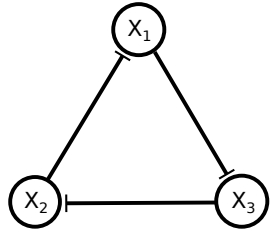}
\end{center}
\caption{A schematic representation of the repressilator circuit.} \label{fig:fig1}
\end{figure}

We first consider the repressilator circuit \cite{Elowitz00}, whose GRN is shown in \cref{fig:fig1}. Here each gene represents a master regulator transcription factor (TFs) of the differentiation process for a specific fate.

Assuming fast (un)binding dynamics of the TFs to DNA, and fast mRNA dynamics, we can reduce the system by adiabatic elimination (see Supplementary Materials) \cite{Greenhill11}, and describe its dynamics with ordinary differential equations (ODEs):
\begin{equation}\label{eq:rep}
\left\lbrace\begin{aligned}
\frac{dx_1}{dt}=b_1 + \frac{g_1}{1+\alpha_1 x_2^{h_1}} - d_1x_1, \\
\frac{dx_2}{dt}=b_2 + \frac{g_2}{1+\alpha_2 x_3^{h_2}} - d_2x_2, \\
\frac{dx_3}{dt}=b_3 + \frac{g_3}{1+\alpha_3 x_1^{h_3}} - d_3x_3.
\end{aligned}\right.
\end{equation}
Here, $x_i$ ($i=1,\,2,\,3$) describe the concentrations of the proteins $\sf X_i$, $b_i$ are background expression rates, $g_i$ are maximal gene expression rates (which combine transcription, translation and mRNA degradation rates), $\alpha_i$ are ratios of binding to unbinding rates (association constants) of TFs to DNA, $h_i$ are Hill coefficients describing possible cooperative effects (like multimerisation of the TFs \cite{Greenhill11}), and $d_i$ are TF degradation rates. 

\begin{figure}[ht]
\begin{center}
\includegraphics{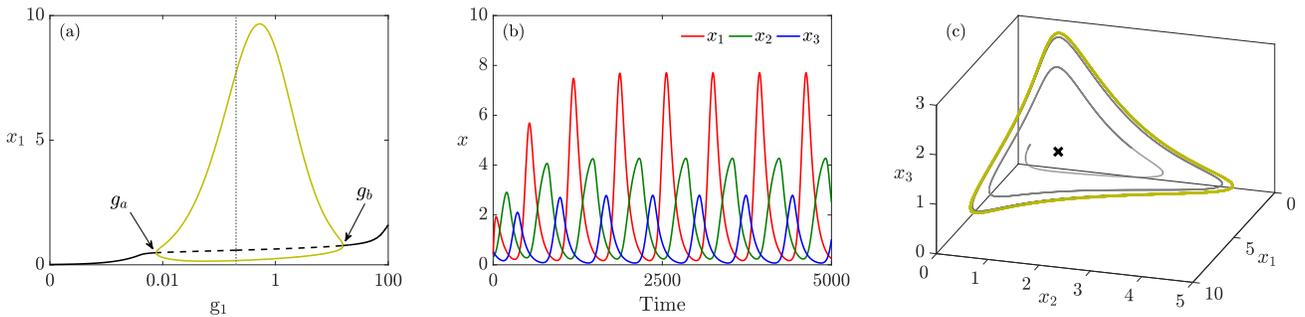}
\end{center}
\caption{a) Bifurcation diagram for the repressilator as $g_1$ varies. A family of stable periodic orbits (olive), shown with maxima and minima of $x_1$, emanates from a supercritical Hopf bifurcation at $g_a$ and terminates at another supercritical Hopf bifurcation at $g_b$. The stable (solid) branch of equilibria (black) becomes unstable (dashed) between these two bifurcations. The dotted vertical indicates the $g_1$-value used for panels (b) and (c). b) Time courses of oscillatory concentrations $x_1$, $x_2$, $x_3$ when $g_1 = 0.2$, $g_2 = g_3 = 0.05$. c) The solution shown in panel (b) in the phase space (grey) converging to the periodic orbit (olive). The cross ($\times$) is the unique saddle equilibrium, at $x_{eq}=(0.58469,\,0.87298,\,0.45578)$. Remaining parameter values are: $b_i = 10^{-5}$, $\alpha_i=50$, $d_i=0.01$, $h_i=3$, for $i=1,2,3$.}
\label{fig:fig2}
\end{figure}

Analysis (see e.g. \cite{Elowitz00}) shows that the repressilator supports oscillatory behaviour over wide regions of parameter space. \cref{fig:fig2} shows (from left) numerical solutions exhibiting the emergence of a family of stable periodic orbits through a Hopf bifurcation as the parameter $g_1$ increases, time series of sustained oscillations, and the limit cycle in phase space. Oscillatory behaviour exists when $h>h^*\approx 2.20498$, for the parameter values used in \cref{fig:fig2}.

For $h=3$, the oscillations exist between the two Hopf bifurcations at $g_a\approx 0.00764$ and $g_b\approx 15.7919$ suggesting that the system is an example of S1 -- variations in a single parameter are sufficient to enter and also to exit the oscillatory regime. Our hypothesis for S1 is that external signalling might increase $g_1=g_1(t)$ over time, letting the system sequentially explore non-oscillatory, oscillatory, and a new non-oscillatory regime. However, exit from the oscillatory regime at $g_1(t) = g_b$ implies convergence onto the single available equilibrium, making the GRN unable to capture the idea of multiple equilibria necessary to describe stem cell differentiation. 

Further, the period of the limit cycles remains bounded as $g_1$ varies. The shape of the limit cycle (\cref{fig:fig2}(c)) shows that the amplitudes of the oscillations for different protein concentrations may differ significantly, and generally appear to spike as the orbit approaches the coordinate axes; two proteins are expressed at relatively low levels while the third is higher. Our simulations show that even when the spiking of gene $\sf X_1$ is more pronounced, with a huge increase in oscillation amplitude, only a logarithimic or algebraic dependency of the oscillation period appears (see Supplementary Materials), demonstrating the absence of any substantial slow-down around the genes maximal expression substates, and failing thereby to effectively describe any of the scenarios S1 or S2. Furthermore,  no mechanism is provided in this model for the stem cell to differentiate into more than one cell type, as no alternative equilibria exist. Thus the standard repressilator does not capture the qualitative changes in behaviour associated with differentiation.

\section{The cross-repressilator with an `OR gate'}
\label{sec:model1}

\begin{figure}[ht]
\begin{center}
\includegraphics{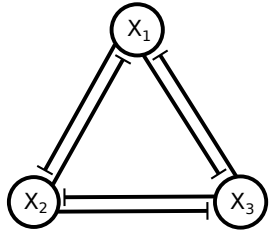}
\end{center}
\caption{Illustration of cross-repression among three genes $\sf X_1$, $\sf X_2$ and $\sf X_3$ formed with two copies of the circuit shown in \cref{fig:fig1}, oriented in opposing directions.} 
\label{fig:fig3}
\end{figure}

We now introduce a variant of the standard repressilator (`cross-repressilator'), where the inhibitory cycle is accompanied by a second inhibitory cycle arranged in the opposite direction, as depicted in \cref{fig:fig3}. We assume that regulation of each gene follows an `OR gate', i.e. each gene is active when \emph{both} inhibitors are absent (see Supplementary Materials).

Making the same assumptions as the standard repressilator, the dynamics of this GRN can be written as
\begin{equation}\label{eq:model1}
\left\lbrace\begin{aligned}
\frac{dx_1}{dt} = b_1 + \frac{g_1}{(1+\alpha_1 x_2^{h_{12}})(1+\beta_1 x_3^{h_{13}})}-d_1x_1, \\
\frac{dx_2}{dt} = b_2 + \frac{g_2}{(1+\alpha_2 x_3^{h_{23}})(1+\beta_2 x_1^{h_{21}})}-d_2x_2, \\
\frac{dx_3}{dt} = b_3 + \frac{g_3}{(1+\alpha_3 x_1^{h_{31}})(1+\beta_3 x_2^{h_{32}})}-d_3x_3,
\end{aligned}\right.
\end{equation}
where the repression of each gene is described by the multiplication of two decreasing Hill functions to capture the assumed 	`OR gate', (see Supplementary Materials). Like system~\cref{eq:rep}, the last term in each equation represents protein degradation. For simplicity, and since it does not cause qualitative changes in the system dynamics, we remove basal expression here and from now throughout this paper. Also, in the numerical simulations we set parameter values equal for each gene, i.e. $g_i=g$, $\alpha_i=\alpha$, $\beta_i=\beta$, $d_i=d$, and $h_{ij}=h$ for $i,j=1,2,3$; we return to this in \cref{sec:break}.

\begin{figure}[ht]
\begin{center}
\includegraphics{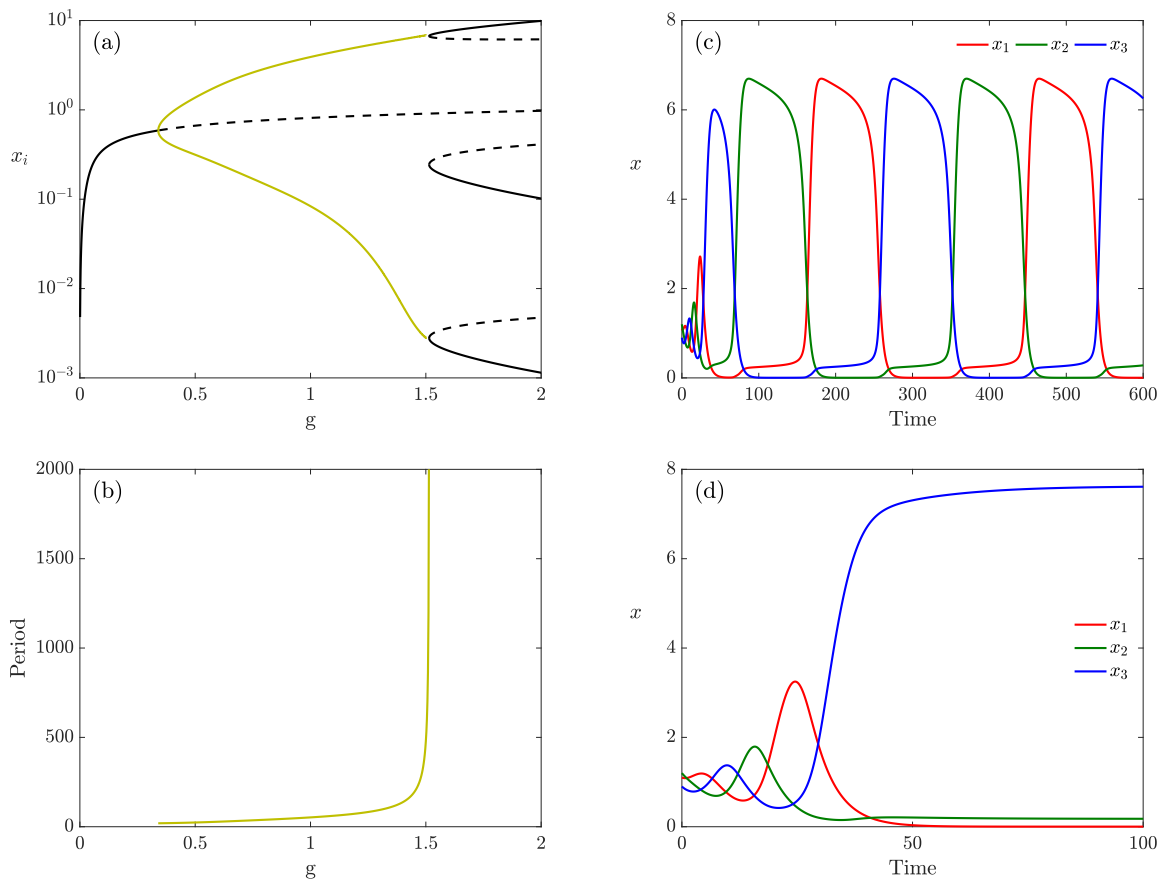}
\end{center}
\caption{a) Bifurcation diagram with respect to $g$ in log-scale for $\alpha=9$, $\beta=0.1$, $h=3$, and $d=0.2$. b) The period of the periodic orbits shown in panel (a) increases to infinity while $g$ tends to the SNIC bifurcation. c) Time courses of the expression level of genes $\sf X_1$ (red), $\sf X_2$ (green) and $\sf X_3$ (blue) showing the sustained oscillations, for $g=1.5$. d) For $g=1.6$, the system stops oscillating and exits to a differentiated cell.}
\label{fig:fig4}
\end{figure}

Numerical simulations of ODEs~\cref{eq:model1} are shown in \cref{fig:fig4}. \Cref{fig:fig4}(a) shows the bifurcation diagram for system~\cref{eq:model1} with respect to the common parameter $g$, all other parameters fixed. The (dashed) solid curves indicate (un)stable equilibria (black) or (maxima and minima of) periodic orbits (olive green). For $g \lessapprox 0.4$, system~\cref{eq:model1} has a single stable equilibrium which is (due to the equalities in parameter values) symmetric: $x_1=x_2=x_3$. The cyclic symmetry apparent in \cref{eq:model1} implies that the Jacobian matrix evaluated at this symmetric equilibrium has a pair of eigenvalues with equal real parts (see Supplementary Materials). When this complex conjugate pair crosses the imaginary axis, a supercritical Hopf bifurcation occurs generating a family of stable periodic orbits. As $g$ increases further, the stable periodic orbits disappear at a global bifurcation involving three new equilibria created at saddle-node bifurcations.

Moreover, there are two distinct global bifurcation mechanisms that remove the periodic orbits. For $\alpha=9$ (shown in \cref{fig:fig4}), we find that the periodic orbits disappear at a SNIC (Saddle-Node on an Invariant Cycle) bifurcation where three saddle-node bifurcations, related by the cyclic symmetry, occur on the periodic orbit simultaneously. The SNIC bifurcation is well-known in simple dynamical systems, such as the nonlinear pendulum with forcing and damping \cite{Coullet2005} and in oscillators under external periodic perturbations \cite{Andronov1996,Tsai2013}. For $g > g_{\rm SNIC}\approx 1.51449$, system~\cref{eq:model1} has no periodic orbit but three equilibria related by the cyclic symmetry. For each of the equilibria, one of the $x_i$ is significantly larger than the other two. Initial conditions (ICs) determine which of these equilibria attract the trajectory.

\Cref{fig:fig4}(b) shows that, in contrast with the standard repressilator, the period of the periodic orbits diverges as $g$ approaches $g_{\rm SNIC}$. Typical trajectories spend longer and longer in the vicinity of the three `slow regions' of phase space where the saddle-node bifurcations are about to occur, and exhibit rapid transitions between these `slow regions'. This phenomenon has been referred to as the `ghost of SNIC' \cite{Strogatz}, and indicates the possibility of S2. We identify these slow regions with the sub-states discussed earlier.

In \cref{fig:fig4}~(right), we show time courses of system~\cref{eq:model1} before~(c) and after~(d) the SNIC bifurcation with same ICs. \Cref{fig:fig4}(c) shows the time evolution of the three genes for g$g = 1.5 < g_{\rm SNIC}$. Whenever one of the genes is highly expressed, the expression levels of the other two is low. In \cref{fig:fig4}(d), we observe the exit from the oscillatory behaviour associated with a further increase of parameter $g$ ($g=1.6$). Here, gene $\sf X_1$ is stably expressed at a high level, while genes $\sf X_2$ and $\sf X_3$ are stably expressed at very low levels.

In summary, for $g$ increasing, the system shows three different behaviours: (i) convergence to a (unique) stable equilibrium; (ii) oscillations among three sub-states; and (iii) attraction to one of the stable equilibria. The natural interpretation, compatible with S1, is to see the first two cases as reflecting multipotency within the stem cell, while the third case represents attainment of one of multiple differentiated cell types. 

\begin{figure}[ht]
\begin{center}
\includegraphics{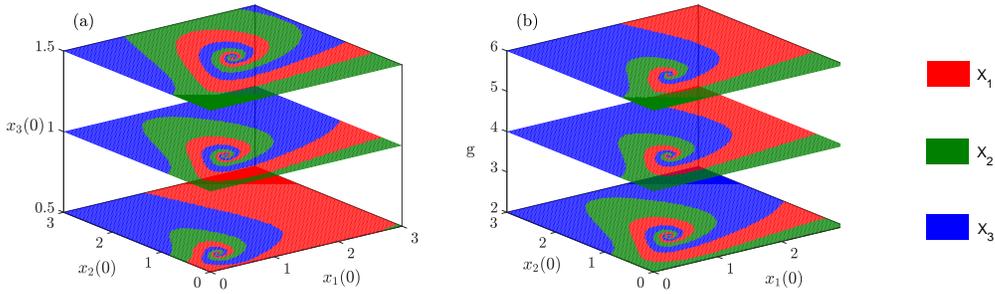}
\end{center}
\caption{a) Colourmaps indicating which gene is expressed starting from ICs in the $(x_1,\,x_2)$-plane, for $g=2$ and three different ICs of $x_3$: $x_3(0)=0.5, 1$ and $1.5$. b) The same colourmaps as panel (a) for $x_3(0)=1$ and three different values of $g$: $g=2,4$ and $6$.}
\label{fig:fig5}
\end{figure}

The phenotype attained in the differentiated state depends on the IC used for the simulation and its intrinsic properties, here illustrated through different $g$-values, as shown in \cref{fig:fig5}. \Cref{fig:fig5}(a) colour-codes the $(x_1,\,x_2)$-plane according to the final equilibrium for three different ICs for $x_3(t)$, and fixed $g$. \Cref{fig:fig5}(b) shows the effect of varying parameter $g$ on the same colourmap when $x_3(0)$ is fixed. The three coloured regions in the $(x_1,\,x_2)$-plane twists around a central point. Increasing $x_3(0)$ tightens the twists of these spiralling regions while the opposite occurs when $g$ increases. Moreover, variations in $g$ keep the central twisting point fixed at $(x_1(0),\,x_2(0))=(1,1)$ whereas increasing $x_3(0)$ makes this point drift away.

\begin{figure}[ht]
\begin{center}
\includegraphics{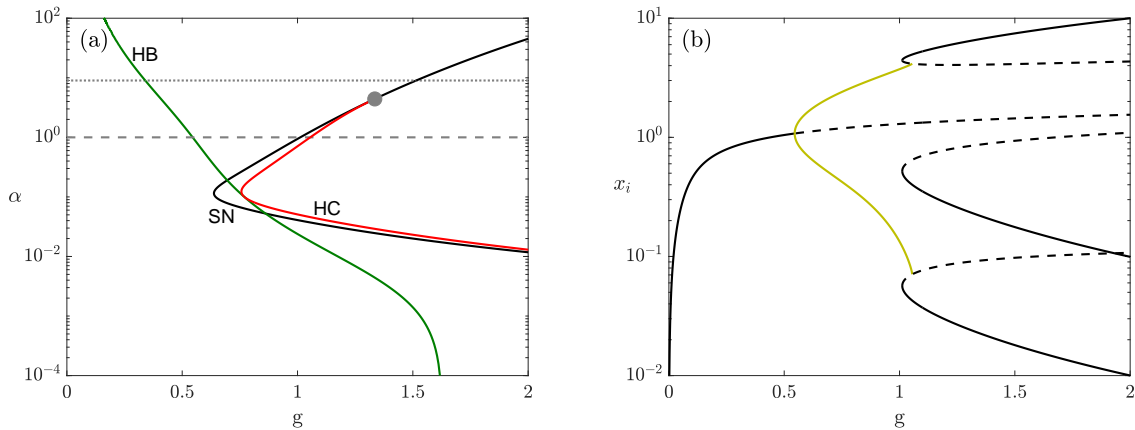}
\end{center}
\caption{a) Two-parameter continuation of saddle-node (black), Hopf (dark green) and heteroclinic (red) bifurcations, shown in \cref{fig:fig4}(a), in the $(g,\alpha)$-plane; vertical axis ($\alpha$) is shown in log-scale. The dotted and dashed lines indicate the corresponding $\alpha$ values used in \cref{fig:fig4}(a) ($\alpha=9$) and in panel (b) ($\alpha=1$), respectively. b) Bifurcation diagram of system~\cref{eq:model1} with respect to $g$ in log-scale when $\alpha=1$. The diagram structure is mainly similar to the one in \cref{fig:fig4}(a). However, here, a small region of multi-stability is created between the oscillatory regime and the stable equilibria such that the family of periodic orbits terminates at a heteroclinic bifurcation (instead of a SNIC).}
\label{fig:fig6}
\end{figure}

To obtain insight into the feasibility of S2 for this system, we explore the bifurcation behaviour in the $(g,\,\alpha)$-plane. \cref{fig:fig6}(a) shows the continuation of the bifurcations identified for $\alpha=9$ in the $(g,\,\alpha)$-plane. When $\alpha$ decreases, the location of the Hopf bifurcation curve ($\sf HB$) moves to larger $g$ whereas the curve of saddle-node bifurcations ($\sf SN$) on which the SNIC bifurcation occurs moves to lower values of g. At $\alpha=\alpha^\ast \approx 3$, the SNIC regime terminates at a codimension-two bifurcation, indicated by a grey dot. Two new bifurcation curves emerge from this point, indicated by the labels $\sf HC$  (homoclinic) and $\sf SN$ (saddle-node). For $\alpha <\alpha^\ast$ the three saddle-node bifurcations (due to symmetry) occur away from the periodic orbit in phase space, resulting in a region of the $(g,\,\alpha)$-plane in which these new equilibria (three stable and three saddle equilibria) coexist with the stable periodic orbits. The periodic orbits then collide with the three saddle equilibria at a global bifurcation and disappear. As in the SNIC bifurcation, the period of the periodic orbits goes to infinity as the orbits approach the global bifurcation, but the mathematical details of the behaviour differ from the SNIC case, as shown for example by the scaling of the period with distance to the bifurcation point being different (see Supplementary Materials).

We observe also that as $\alpha$ decreases further the heteroclinic bifurcation ($\sf HC$) and $\sf SN$ curves move further apart from each other in the $(g,\alpha)$-plane, with maximum separation at $\alpha \approx 0.1$.

For $\alpha<\alpha^\ast$ the $\sf HB$ (green) and $\sf HC$ curves (red) determine the left and right boundaries of the oscillatory region, respectively. Also note that the $\sf HB$ and $\sf HC$ curves meet and are tangent to each other at $\alpha=0.1$ which is not generic but arises here due to the additional symmetry existing when $\alpha=0.1=\beta$. When $\alpha=\beta$, ODEs~\cref{eq:model1} are symmetric under two independent symmetry operations: the cyclic permutation symmetry $(x_1,x_2,x_3) \rightarrow (x_2,x_3,x_1)$ and the interchange of any two of the three variables, e.g. $(x_1,x_2,x_3) \rightarrow (x_2,x_1,x_3)$. In group theory notation, the symmetry group is now the symmetries of an equilateral triangle $D_3$ rather than just the cyclic group of order three, $\mathbb{Z}_3$. This additional symmetry implies the tangency between the two curves. The analysis around this point will be discussed in detail in a separate paper \cite{Dawes2021} and has been observed in other contexts \cite{gambaudo1985,marques2013}.

\Cref{fig:fig6}(b) shows the bifurcation diagram, varying $g$, for $\alpha=1$ corresponding to the lower horizontal dashed line indicated in \cref{fig:fig6}(a). Here we observe an interval of $g$-values in which the stable periodic orbit coexists with the three stable equilibria. Consequently we would expect that this transition to a stable equilibrium occurs rather smoothly since the sub-states already exist prior to the termination of the oscillations. 

These observations support S2. Once the system has entered the oscillatory regime by tuning of $g$ through a first signalling pathway, a second signalling pathway might intervene to lower $\alpha$, stop the oscillations, and drive the system to differentiation. This scenario requires that the change of $\alpha$ happens faster than the permanence of the system close to a selected sub-state. This slowdown of the cycling dynamics is therefore essential for this mechanism. 

Thus, after establishment of the stem cell, increasing the first signal initiates the cyclical expression of fate- specific TFs corresponding to the different differentiated cell-types. Subsequently, in each primed sub-state, as a result of TF specific activation of sensitivity to fate specification signals, the cell is sensitised to local signals. When a cell receives such a signal for sufficient time, this shifts $\alpha$ and removes the oscillations, thereby initiating differentiation.

\subsection{Illustrating the two scenarios with time-dependent parameters }
\label{sec:time}

\begin{figure}[ht]
\begin{center}
\includegraphics{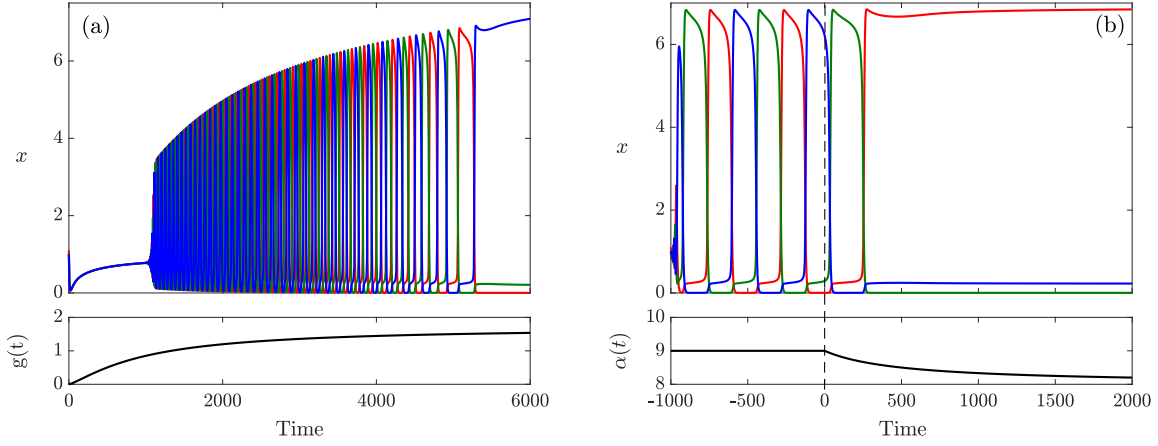}
\end{center}
\caption{a) The development over time from a multipotent cell to a differentiated cell controlled by a time-dependent profile $g(t)$ (shown in the lower panel) according to S1. Parameters are as \cref{fig:fig4}(c). b) Development from a stable oscillatory regime into a differentiated cell according to S2, when the secondary parameter $\alpha$ starts decreasing at $\tilde{t}=0$ (shown in the lower panel). Here $g = 1.5$ is held constant, and all other parameters take the same values as in \cref{fig:fig4}(c).}
\label{fig:fig7}
\end{figure}

So far we have assumed that all the parameters of system~\cref{eq:model1} are fixed and do not vary in time. In this section, we explore more directly the two scenarios S1 and S2 by assuming that relevant parameters are regulated in time by an external signalling pathway. We will consider the case when either one or both of $g$ and $\alpha$ are time dependent. In S1, parameter $g$ only is responsible for shifting the cell from the multipotent regime into the oscillatory regime and then into the differentiated cell. Accordingly, we assume that $g(t)$ follows the form of an increasing Hill function, from 0 and saturating at $\hat{g}>g_{\rm SNIC}$ as
$$g(t) = \frac{\hat{g}t^p}{\tau^p+t^p},$$
where we set $\hat{g} = 1.7$, $p = 1.25$ and $\tau = 1000$. The bottom panel of \cref{fig:fig7}(a) illustrates the time-variation of $g(t)$. The top panel in \cref{fig:fig7}(a) shows how the system responds to $g(t)$, moving through the three different regimes as $g(t)$ increases and tends to the value $\hat{g}$ (compare with \cref{fig:fig4}(a)). In detail, the system first settles quickly from the IC to the single stable equilibrium and follows it as $g(t)$ slowly approaches the Hopf bifurcation. After Hopf we observe the onset of rapid but finite-frequency oscillations. As $g\longrightarrow g_{\rm SNIC}$, the period of oscillations increases significantly, giving rise to dynamical phases close to each of the three sub-states. Finally, when $g(t)$ crosses $g_{\rm SNIC}$, the oscillations disappear and the system settles at one of the three equilibria, each corresponding to a different expressed gene. For example, in \cref{fig:fig7}, gene $\sf X_3$ (blue) is expressed and the other two genes remain at a low level.

As discussed for S1, the precise state of the system when $g$ goes through  $g_{\rm SNIC}$, as well as the value of the saturation level $\hat{g}$ play crucial roles in selecting cell fate. The shape of the profile of $g(t)$ is also relevant for this selection.

For S2, we hypothesise that two parameters, $g$ and $\alpha$, are involved in cell differentiation. As shown in \cref{fig:fig7}(b), we first assume that parameter $g$ shifts the system into the oscillatory regime while $\alpha$ remains fixed. At a certain time, we suppose a second pathway causes a decrease in $\alpha$, modelled as
$$\alpha(t) = 9 - \mathcal{H}(t-\tilde{t}) \times \frac{(t-\tilde{t})^p}{\tau^p+(t-\tilde{t})^p},$$
where $\mathcal{H}$ is a Heaviside step function (i.e. $\mathcal{H}(s)=0$ when $s<0$ and $\mathcal{H}(s)=1$ when $s>0$), $\tilde{t}$ is the time when the time dependency of $\alpha(t)$ starts, and the constants $\tau = 1000$ and $p = 1$ control the profile shape.

\begin{figure}[ht]
\begin{center}
\includegraphics{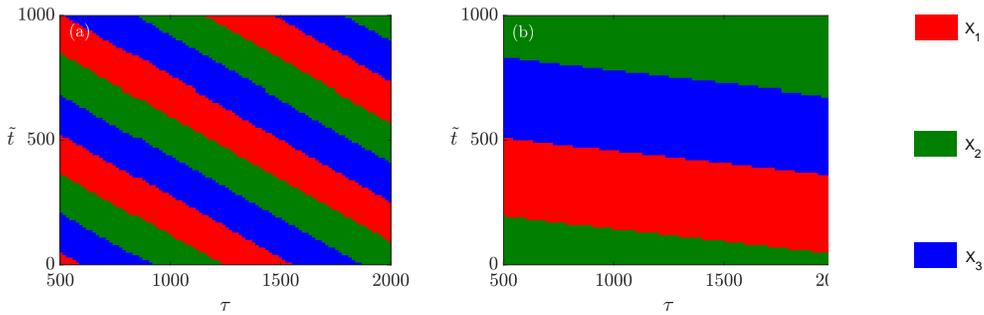}
\end{center}
\caption{Colourmaps show the expressed gene in $(\tau,\tilde{t})$-plane starting from the same IC for $g = 1.5$~(a) and $g = 1.51$~(b).}
\label{fig:fig8}
\end{figure}

S2 clearly allows more control over the final fate of the cell through the increased number of coefficients involved in specifying the time-dependency of the second bifurcation parameter. To understand these further dependencies, we analyse the effect of variations in the switch-on time $\tilde{t}$ and the time scale parameter $\tau$ for the saturation of $\alpha(t)$. Shown in \cref{fig:fig8}, we colour-code the parameter plane $(\tau,\tilde{t})$ based on the expressed gene for $g = 1.5$ (a) and $g = 1.51$ (b). The $(\tau,\tilde{t})$-plane is divided into stripes which periodically alternate between red, green and blue. The lower value of $g$ in \cref{fig:fig8}(a) corresponds to narrower stripes as a result of the smaller oscillation period, with the system being more sensitive to changes in switch-on time and saturation time scale. Further, more horizontal stripes in \cref{fig:fig8}(b) indicate that cell fates become less sensitive to $\tau$ for this larger value of $g$, at least in the regime close to the global bifurcation.

\subsection{Breaking the symmetry}
\label{sec:break}

\begin{figure}[ht]
\begin{center}
\includegraphics{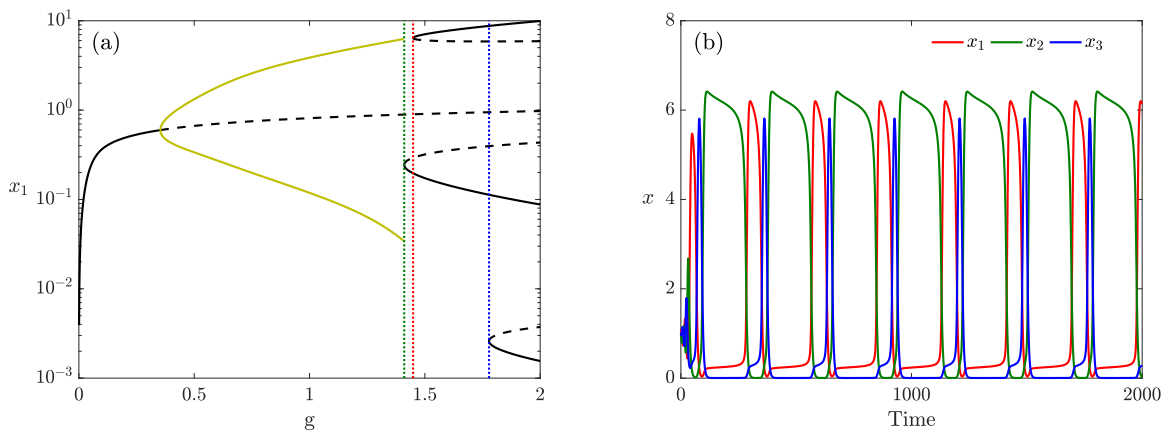}
\end{center}
\caption{a) Bifurcation diagram of the system for $\alpha_1=7$, $\alpha_2=9$, $\alpha_3=8$, $\beta_1=0.115$, $\beta_2=0.07$, $\beta_3=0.1$ and $d=0.2$ in log-scale. b) Time courses of the oscillations in the expression level of genes $\sf X_1$ (red), $\sf X_2$ (green) and $\sf X_3$ (blue) for the same parameters in panel (a) and $g=1.4$.}
\label{fig:fig9}
\end{figure}

Until now, we assumed that system~\cref{eq:model1} is symmetric and all the association rates $\alpha_i$ and $\beta_j$, $i,j=1,2,3$ for the TFs $\sf X_1$, $\sf X_2$ and $\sf X_3$ are equal. Here, we break the symmetry in the reciprocal repressing loops in \cref{fig:fig3} and consider the case that $\alpha_i$ and $\beta_j$ take different values, as specified in \cref{fig:fig9}. Note that we keep the remaining parameters $g$ and $d$, unchanged.

\Cref{fig:fig9}(a) shows the bifurcation diagram with respect to $g$ in log-scale. The structure remains qualitatively unchanged compared to the symmetric case. For $g$ small, a unique stable equilibrium exists, which undergoes a Hopf bifurcation and becomes unstable. For $g$ large, there are six branches of equilibria; the stable branches merge with the saddle-type ones and disappear as $g$ decreases. The family of stable periodic orbits emerging from the Hopf bifurcation terminates at a SNIC bifurcation. However, here, only one equilibrium lies on the invariant cycle contrary to the symmetric case, with three equilibria on the invariant cycle. This happens because the saddle-node bifurcations, at which the six branches of equilibria appear, occur at different $g$-values (they are no longer symmetry-related), indicated by the vertical dotted lines.

\Cref{fig:fig9}(b) shows time courses of the genes concentrations just before the SNIC bifurcation at $g_{\rm SNIC}\approx 1.41054$; the other parameter values are as in \cref{fig:fig9}(a). As before, each TF is expressed cyclically when the other two are at low levels. However, gene $\sf X_2$ remains expressed for a significantly longer time than the other two, particularly $\sf X_3$. Due to the occurrence of the saddle-node bifurcation with high level of $\sf X_2$ at smaller values of $g$, the slow region created around the bifurcation has a stronger effect which leads to a longer `dwelling time' at expressed $\sf X_2$. In fact, the dwelling time for each expressed gene depends on the distance of g-value at which the corresponding saddle-node bifurcation (indicated with a dotted vertical line of the same colour) occurs. 

Analysis of \cref{fig:fig9} determines that with the current parameters, the conditions favour gene $\sf X_2$; therefore, when $g$ crosses to the right of the SNIC bifurcation, the system settles at a stable equilibrium where $\sf X_2$ is expressed. Changing this condition results in a shift in the location of the saddle-node bifurcations along the horizontal axis and might alter the order of saddle-node bifurcations; consequently, lengthening the dwelling time for another gene.

\section{The cross-repressilator with an `AND gate'}
\label{sec:model2}

The third model we consider here is similar to system~\cref{eq:model1} with the difference that each gene is repressed only if all other genes are highly expressed (a so-called `AND gate').

This model is described by
\begin{equation}\label{eq:model2}
\left\lbrace\begin{aligned}
\frac{dx_1}{dt} = b_1 +  \frac{g_{11}+g_{12}x_2^{h_{12}}+g_{13}x_3^{h_{13}}}{(1+\alpha_1x_2^{h_{12}})(1+\beta_1x_3^{h_{13}})}-d_1x_1, \\
\frac{dx_2}{dt} = b_2 + \frac{g_{21}+g_{22}x_3^{h_{23}}+g_{23}x_1^{h_{21}}}{(1+\alpha_2x_3^{h_{23}})(1+\beta_2x_1^{h_{21}})}-d_2x_2, \\
\frac{dx_3}{dt} = b_3 + \frac{g_{31}+g_{32}x_1^{h_{31}}+g_{33}x_2^{h_{32}}}{(1+\alpha_3x_1^{h_{31}})(1+\beta_3x_2^{h_{32}})}-d_3x_3,
\end{aligned}\right.
\end{equation}
where again the decreasing Hill functions describe the effect of repressing genes (see Supplementary Materials). We also set $g_{ij}=g_j$ for $i,j=1,\ldots,3$ and $d_1=d_2=d_3=d$. We also fix $\alpha_i=\beta_i=1$ and the Hill function exponents $h_{ij}=3$ for $i,j=1,2,3$. 

\begin{figure}[ht]
\begin{center}
\includegraphics{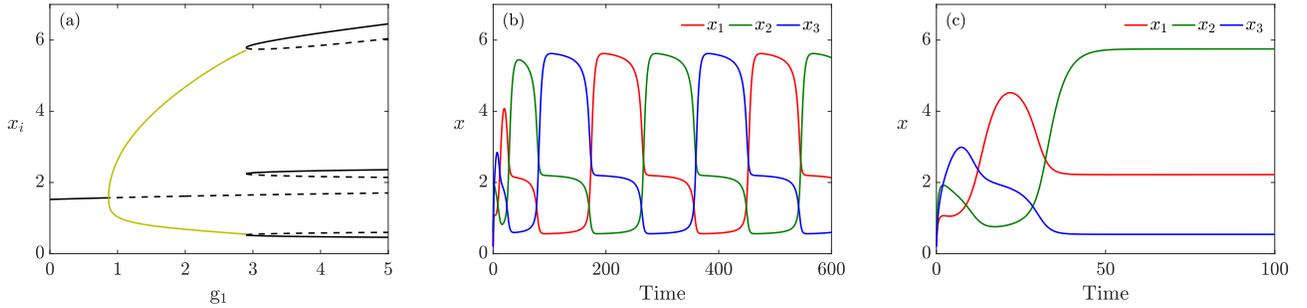}
\end{center}
\caption{: a) Bifurcation diagram of system~\cref{eq:model3} with respect to $g_1$ for $g_2=0.76$, $g_3=1.9$ and $d=0.3$. b) Time courses of the expression level of genes $\sf X_1$ (red), $\sf X_2$ (green) and $\sf X_3$ (blue) showing the emergence of an intermediate expression level in the oscillatory regime of system~\cref{eq:model2} when $g_1=2.8$. c) Gene expression levels converging to equilibrium in the tristable regime for $g_1=3$.}
\label{fig:fig10}
\end{figure}

As previously, we explore system~\cref{eq:model2} numerically and compute the bifurcation diagram, varying $g_1$ while keeping $g_2$ and $g_3$ fixed (and unequal). \Cref{fig:fig10}~(a) shows a bifurcation diagram for system~\cref{eq:model2} with respect to $g_1$. For small values of $g_1$, the system has a single stable equilibrium (black solid line) which undergoes a supercritical Hopf bifurcation and becomes unstable (dashed line). The stable periodic orbits (olive) are indicated by the maximum and minimum values of $x_1$. The periodic orbit, similar to the `OR gate' case, terminates at a SNIC bifurcation where it collides with a saddle-node bifurcation. Note that these branches are part of a single branch of equilibria terminating at another oscillatory regime at very large parameter values $g_1$ (not shown). Panels~(b) and~(c) in \cref{fig:fig10} show time courses of the expression level of each gene in system~\cref{eq:model2} for values of $g_1$ mentioned in the caption.

Clearly, the behaviour of the `AND gate' model follows closely the analysis presented for system~\cref{eq:model1}, but with one important and biologically relevant difference. Compared to system~\cref{eq:model1}, in which only one gene could be highly expressed and the other two genes were necessarily at very low levels, in system (3), two genes can be expressed at high or moderate levels simultaneously. Of course, the precise levels depend on parameter values, but the `AND gate' structure is crucial in admitting states with more than one high expression level. The form of the periodic orbit also changes and shows that it lingers in sub-states where two gene expression levels are non-zero. This notable feature is potentially useful in matching models with quantitative experimental data where it allows one to select a model that allows or prevents the existence of phenotypes dependent upon the expression of a single or a combination of master regulator TFs.

\section{Extension to four and five genes }
\label{sec:four_genes}

\begin{figure}[ht]
\begin{center}
\includegraphics{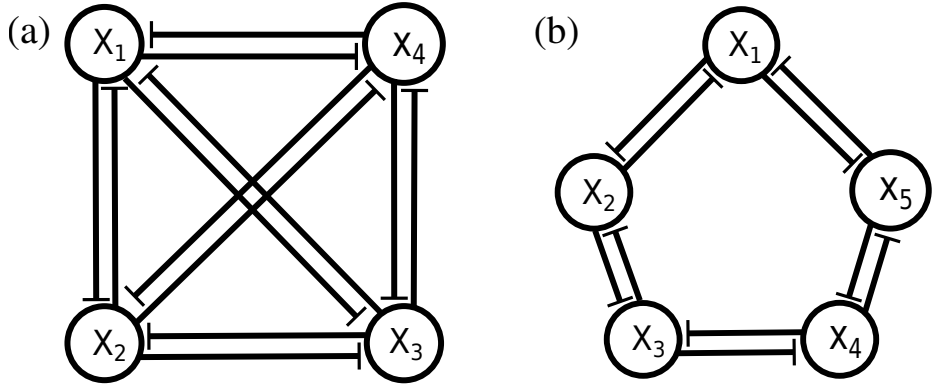}
\end{center}
\caption{a) Extension of the cross-repressilator GRN to four genes with cross connections that support oscillations. b) A cyclical arrangement of five genes, sufficient to drive oscillations without additional cross connections.}
\label{fig:fig11}
\end{figure}

So far, we considered GRNs with only three genes. We now investigate the dynamics with more than three genes. For the same cross-repressilator in `AND gate' configuration discussed in the previous section, an odd number of genes is necessary to create a negative feedback loop that supports sustained oscillations \cite{Thomas81}. Hence a trivial extension of the cross-repressilator to four genes, preserving the cyclical topology, fails to reproduce oscillatory dynamics, and therefore we do not consider it here. Even numbers of genes in the cycle result in states where the gene expression levels alternate around the loop. However, adding cross connections to this four-gene network so as to build three-gene negative loops, embedded within the four-gene GRN as shown in \cref{fig:fig11}(a), recovers the oscillatory dynamics. For five genes, the simple cyclical topology of the cross-repressilator, with only `nearest neighbour' interactions and without additional cross-regulation, shown in \cref{fig:fig11}(b), maintains oscillatory behaviour.

\subsection{The four gene network}

Here, we describe the dynamics of a GRN with four genes, as illustrated in \cref{fig:fig11}(a). The network is fully connected such that each gene is repressed by the other three genes, operating in `AND gate' configuration. The results in (see Supplementary Materials for details):
\begin{equation}\label{eq:model3}
\left\lbrace\begin{aligned}
\frac{dx_1}{dt} = b_1 +  \frac{g_{11}+g_{12}x_2^{h_{12}}+g_{13}x_3^{h_{13}}+g_{14}x_4^{h_{14}}+g_{15}x_2^{h_{12}}x_3^{h_{13}}+g_{16}x_2^{h_{12}}x_4^{h_{14}}+g_{17}x_3^{h_{13}}x_4^{h_{14}}}{(1+\alpha_1x_2^{h_{12}})(1+\beta_1x_3^{h_{13}})(1+\gamma_1x_4^{h_{14}})}-d_1x_2, \\
\frac{dx_2}{dt} = b_2 + \frac{g_{21}+g_{22}x_3^{h_{23}}+g_{23}x_4^{h_{24}}+g_{24}x_1^{h_{21}}+g_{25}x_3^{h_{23}}x_4^{h_{24}}+g_{26}x_3^{h_{23}}x_1^{h_{21}}+g_{27}x_4^{h_{24}}x_1^{h_{21}}}{(1+\alpha_2x_3^{h_{23}})(1+\beta_2x_4^{h_{24}})(1+\gamma_2x_1^{h_{21}})}-d_2x_2, \\
\frac{dx_3}{dt} = b_3 + \frac{g_{31}+g_{32}x_4^{h_{34}}+g_{33}x_1^{h_{31}}+g_{34}x_2^{h_{32}}+g_{35}x_4^{h_{34}}x_1^{h_{31}}+g_{36}x_4^{h_{34}}x_2^{h_{32}}+g_{37}x_1^{h_{31}}x_2^{h_{32}}}{(1+\alpha_3x_4^{h_{34}})(1+\beta_3x_1^{h_{31}})(1+\gamma_3x_2^{h_{32}})}-d_3x_3, \\
\frac{dx_4}{dt} = b_4 +  \frac{g_{41}+g_{42}x_1^{h_{41}}+g_{43}x_2^{h_{42}}+g_{44}x_3^{h_{43}}+g_{45}x_1^{h_{41}}x_2^{h_{42}}+g_{46}x_1^{h_{41}}x_3^{h_{43}}+g_{47}x_2^{h_{42}}x_3^{h_{43}}}{(1+\alpha_4x_1^{h_{41}})(1+\beta_4x_2^{h_{42}})(1+\gamma_4x_3^{h_{43}})}-d_4x_4.
\end{aligned}\right.
\end{equation}

For simplicity, we again set $g_{ij}=g_j$ and $d_i = d$ for $i,j=1,\ldots,4$. We also keep fixed the values $g_j=0.5$ for $j=5,\ldots,7$, $\alpha_i=\beta_i=\gamma_i=1$ and $h_{ij}=3$ for $i,j=1,\ldots,4$ ($i \neq j$). 

\begin{figure}[ht]
\begin{center}
\includegraphics{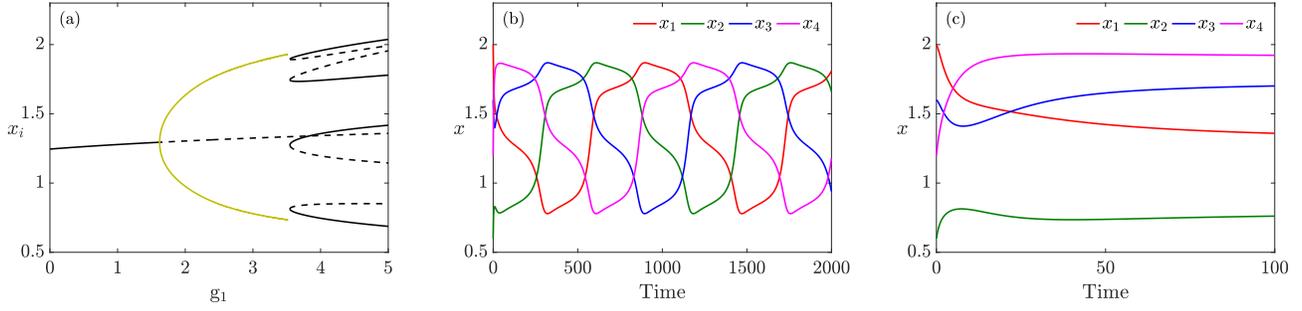}
\end{center}
\caption{a) Bifurcation diagram of system~\cref{eq:model3} with respect to $g_1$ for $g_2=1.6$, $g_3=0.5$, $g_4=1.5$ and $d=0.4$. (b) and (c) Time courses of the expression level of genes $\sf X_1$ (red), $\sf X_2$ (green), $\sf X_3$ (blue) and $\sf X_4$ (magenta) with two different behaviours: in the multipotent oscillatory regime of system~\cref{eq:model3} when $g_1=3$~(b) and in the multistable regime for $g_1=3.6$~(c).}
\label{fig:fig12}
\end{figure}

\Cref{fig:fig12}(a) shows the bifurcation diagram of system~\cref{eq:model3} with respect to $g_1$, together with time courses of the expression levels of all genes in the system for two specific values of parameter $g_1$ in the multipotent~(b) and differentiated~(c) regimes. In \cref{fig:fig12}(a), a single branch of stable equilibria (solid; black) undergoes a Hopf bifurcation and becomes unstable/saddle (dashed). The envelope of stable periodic orbits (olive) emanating from the Hopf bifurcation terminates at a SNIC bifurcation where new branches of equilibria appear. Four of these new branches are stable, and collide with the other four unstable branches, disappearing at saddle-node bifurcations.

Panels~(b) and~(c) in \cref{fig:fig12} show time series of the expression levels of the genes in the oscillatory and differentiated regimes for $g_1=3$ and $g_1=3.6$, respectively. One notable, and biologically relevant, feature of this model is that all four genes can be expressed, with one gene highly expressed, two genes expressed at intermediate levels, and one gene only slightly expressed. This scenario is consistent with differentiation reflecting a combinatorial role for multiple TFs. The relative concentrations of the expressed genes would open up possibilities for differentiation to multiple cell types.

\subsection{The five gene network}

\begin{figure}[ht]
\begin{center}
\includegraphics{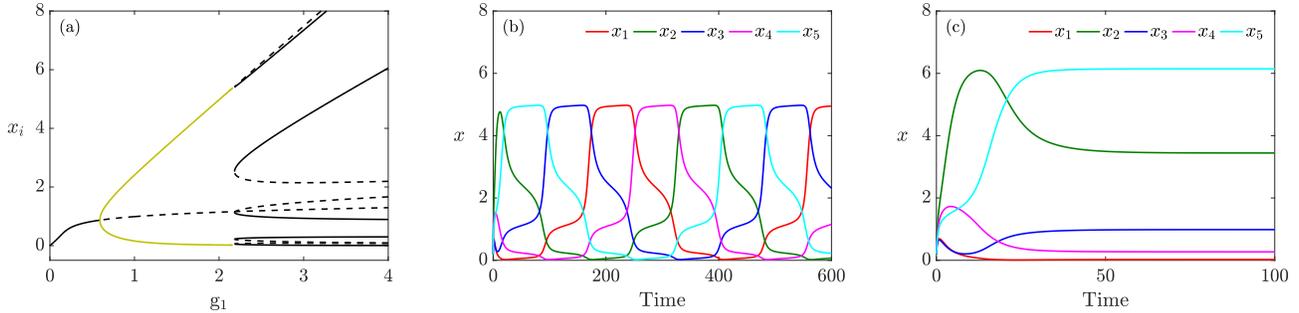}
\end{center}
\caption{(a) Bifurcation diagram of system (5) with respect to $g_1$ when $g_2=0.35$, $g_3=0.2$ and $d=0.4$. (b) Time series of the solutions in the oscillatory/multipotent regime when $g_1=2$. (c), and multistable/ differentiated regime for $g_1=2.5$.}
\label{fig:fig13}
\end{figure}

We now consider the five-gene network depicted in \cref{fig:fig11}(b). The network comprises two inhibitory loops arranged in opposing directions so that gene $X_i$ is repressed by genes $X_{i-1}$ and $X_{i+1}$, again in `AND gate' configuration. Therefore, similar to the three-gene case, the dynamical equations are
\begin{equation}\label{eq:model5}
\left\lbrace\begin{aligned}
\frac{dx_1}{dt} = b_1 +  \frac{g_{11}+g_{12}x_2^{h_{12}}+g_{13}x_5^{h_{15}}}{(1+\alpha_1x_2^{h_{12}})(1+\beta_1x_5^{h_{15}})}-d_1x_1, \\
\frac{dx_2}{dt} = b_2 +  \frac{g_{21}+g_{22}x_3^{h_{23}}+g_{23}x_1^{h_{21}}}{(1+\alpha_2x_3^{h_{23}})(1+\beta_2x_1^{h_{21}})}-d_2x_2, \\
\frac{dx_3}{dt} = b_3 +  \frac{g_{31}+g_{32}x_4^{h_{31}}+g_{33}x_2^{h_{32}}}{(1+\alpha_3x_4^{h_{34}})(1+\beta_3x_2^{h_{32}})}-d_3x_3, \\
\frac{dx_4}{dt} = b_4 +  \frac{g_{41}+g_{42}x_5^{h_{43}}+g_{43}x_3^{h_{43}}}{(1+\alpha_4x_5^{h_{45}})(1+\beta_4x_3^{h_{43}})}-d_4x_4, \\
\frac{dx_5}{dt} = b_5 +  \frac{g_{51}+g_{52}x_1^{h_{51}}+g_{53}x_4^{h_{54}}}{(1+\alpha_5x_1^{h_{51}})(1+\beta_5x_4^{h_{54}})}-d_5x_5.
\end{aligned}\right.
\end{equation}
where the maximal expression rates imposed by the anticlockwise and the clockwise loops respectively equal to $g_1$, $g_2$ and $d_i$ to $d$. We also fix the values $\alpha_i=\beta_i=\gamma_i=1$ and $h_{ij}=3$ for $i,j=1,\ldots,5$ ($i \neq j$). \Cref{fig:fig13}(a) shows the bifurcation diagram remains unchanged with respect to $g_1$ except for the number of stable (solid; black) and saddle (dashed; black) branches in the multistable regime. \Cref{fig:fig13}(b) exhibits time courses of gene expression levels when $g_1=2$. \Cref{fig:fig13}(c) shows gene concentrations for $g_1=2.5$ when the cell differentiates. In the differentiated regime, the number of expressed genes increases to three compared with the three-gene case. We conjecture that more generally for a cyclic network of $n$ genes, $\frac{n+1}{2}$ genes are expressed at high or moderate (but unequal) levels and $\frac{n-1}{2}$ are either not expressed at all or are only weakly expressed.

Biologically, the stable generation of cells expressing mixtures of key fate specification TFs shown by these four- and five-gene models matches closely the concept of combinatorial TF expression defining similar cell-types. This is most keenly demonstrated in the peripheral and central nervous systems where distinct neuron types express mixtures of TFs that together specify their individual fate \cite{Lloret-Fernandez18}. It may be that the fully differentiated phenotype represents the impact of a mix of TFs expressed at different levels, even where we tend to think of a single dominant master regulator, such as Mitf in melanocyte development.

\section{Conclusions}
\label{sec:conc}

We have presented a detailed bifurcation analysis of four variants of the standard repressilator model, and discussed their suitability as novel generic models for the process of cellular differentiation. While the original repressilator fails to provide convincing dynamical behaviour, we find that the variants proposed here exhibit biologically relevant features. Our modelling moves beyond Waddington's vision of the epigenetic landscape since such a picture cannot capture highly dynamical internal states of cells which our model places centre stage.

Bifurcation analysis of the systems using \cite{Doedel10} reveals that the oscillatory behaviour emerges from a Hopf bifurcation. At the late stages of multipotency, the oscillating cell visits a set of sub-states where the relevant fate-specific transcription factor(s) are highly expressed relative to the other gene(s). The oscillations terminate at a global bifurcation and the cell settles at a state resembling that of a committed cell-type. From a biological perspective, we can think of the individual fate-specific transcription factors (TFs) as defining individual cell-fates. Where they oscillate, the stem cell retains multipotency (it sequentially transits nearby all states), but at the same time it becomes periodically biased to adopt individual fates due, we propose, to changing exposure to fate-specification signals. When such a signal is received at a sufficiently high level, it may force oscillations to cease and a single TF becomes stably expressed at the expense of the others –-- cell differentiation has begun.

We have focused on two possible scenarios for driving the oscillatory behaviour. In the first scenario, the same pathway is responsible for initiating and terminating the oscillations, by the increase of the chosen parameter. In the second scenario, a first pathway brings the system in the oscillatory region, and a second one intervenes to drive the exit from it. Our results show that the second scenario, however, gives more control over cell fate. Moreover, the timing of controlling signals is essential in our framework, since this combines with ICs to determine the specific TF that remains stably expressed, and thus the differentiated cell-type.

Depending on the type of gates assumed in the network, one or two genes (or more in networks of higher dimensionality) can be co-expressed during both the oscillatory phase and in the selected differentiated cell type. During the oscillations, combinations of genes may be expressed transiently –-- these might correspond to partially-restricted cell intermediates (expressing a subset of fate-specific TFs); however, we emphasise that in our scenarios these retain oscillations and hence must be thought of as fully multipotent, despite their `snapshot appearance'. In the selected differentiated cell-type, the combination and the levels of expressed TFs vary with ICs and timing of the fate-specification factor; thus different cell-types, with specific quantitatively distinguishable combinations of TFs, are generated. Such a model now provides a view of how multipotency might be exhibited in stem cells, in contrast to the series of bipotential intermediates that underpin standard thinking in development. Note that the external signals drive a fully multipotent cell to a specific committed state directly, without further intermediates with restricted potency. Finally, the case with which stable (committed) states consisting of mixed expression of two or more TFs can be generated may well relate to the mechanism producing related cell-types, e.g. neuron types in both the central and peripheral nervous system. 

There is a wider context for the results presented here; the identification of dynamical scenarios typical for repressilator-type GRNs, particularly in view of their proposed role in clock circuits underlying global patterning processes in development \cite{JutrasDube2020}. Although the developmental context considered there is different, the formulation and behaviour of the mathematical models has some points of similarity with our analysis, such as the occurrence of Hopf and SNIC bifurcations, which are generic bifurcations involving time-periodic oscillations.

A key difference with our work is ordering of the different regimes. Model~1 of \cite{JutrasDube2020} shows a progression from oscillations to a single stable equilibrium, and then a saddle-node bifurcation that creates multiple equilibria. Model~2 of \cite{JutrasDube2020} has no Hopf bifurcation but instead a SNIC separating oscillatory from multistable regime; there is no single equilibrium. In contrast, our model robustly indicates a progression from a single equilibrium to oscillations via a Hopf bifurcation, with the oscillations then terminating at a global bifurcation where multiple equilibria emerge. We note that the distinction between the SNIC mechanism and the combination of the global and saddle-node bifurcations is likely to be extremely difficult to detect experimentally due to the high level of resolution required in time, and the effects of noise. It is however mathematically important in terms of comparisons between models.

In terms of model construction a key point of difference is that both Model~1 and Model~2 in \cite{JutrasDube2020} are constructed by interpolating, rather artificially, between two mechanisms potentially driven by distinct and unrelated morphogens. Our dynamical model is structurally more robust and parsimonious as our work has shown that the transition from equilibrium to oscillations whose frequency drops to zero, to `differentiated' equilibria, is a natural result of a single GRN mechanism. Naturally, spatial patterning is not part of our model since we focus on the dynamics within a single cell, but establishing connections between temporal dynamics and spatial patterning is an intriguing research direction that we will explore in a future paper.

\section*{Acknowledgments}
We gratefully acknowledge Laura Hattam and Michael Thomas for helpful preliminary modelling work. RNK gratefully acknowledges support from the Institute for Mathematical Innovation at the University of Bath that helped to advance this project in its early stages. Funding was provided by BBSRC grants BB/S01604X/1 (SF and AR) and BB/S015906/1 (KCS, JHPD, RNK).

\end{document}


\title{Novel Generic Oscillatory Mechanisms in Models for Differentiating Stem Cells\\
{\large Supplementary Material}}

\author{Saeed Farjami, Karen Camargo Sosa, \\ Jonathan H.P. Dawes, Robert N. Kelsh, Andrea Rocco}

\date{}

\maketitle

\section{Model Derivation}
\label{sec:appendix}

Following \cite{Greenhill11}, we here provide a derivation of the Hill functions dynamics utilised in the construction of the models in this paper.

\subsection{Regulation by repressive transcription factors}

Many processes are involved in the regulation of a protein. Here, we consider only few of them and assume that the other processes such as transcription and mRNA translation are implicit. The model includes (un)binding of the transcription factor (TF) $R$, degradation of the protein $P$, and synthesis of $P$ when gene $G$ is in the active state, denoted $G^\ast$, as follows.
\be\label{proc}
R + e_{R} \;\autorightleftharpoons{$\alpha_0$}{$\alpha_1$} \;r_{R},
\hspace{2cm}
G^* \;\autorightarrow{$g$}{}\; P,
\hspace{2cm}
P\; \autorightarrow{$d_P$}{}\; \emptyset,
\ee
where $e_R$ and $r_R$ represent the `empty' and `regulated' regulatory elements of $G$, respectively. The parameters over/under the arrows indicate the rate of the processes.

If $[\cdot]$ denote the number of a chemical per cell, the processes in (\ref{proc}) are described with
\begin{equation}\label{eqproc}
\left\lbrace
\begin{aligned}
\frac{d[r_R]}{dt} = \alpha_0[R](1-[r_R])-\alpha_1[r_R], \\
\frac{d[P]}{dt} = g[e_R]-d_P[P].\phantom{-lllllllllll}
\end{aligned}
\right.
\end{equation}
where $[e_R]+[r_R] = 1$ and $[G^\ast] = [e_R]$.

Binding and unbinding of $R$ to the regulatory
element is significantly faster than the protein regulation and
degradation ($\alpha_0,\,\alpha_1 \gg g,\,d_P$). By quasi-steady-state approximation $d[r_R]/dt=0$, we have
\be\label{final}
\frac{d[P]}{dt} = g\frac{1}{1 + \alpha [R]} - d_P [P], 
\ee
where $\alpha = \alpha_o/\alpha_1$.

Here, the TF $R$ regulates the protein $P$ as a monomer giving rise to a Hill coefficient one. Hill functions with higher coefficients can be  derived when regulation proceeds through binding/unbinding of complexes (dimers, trimers, etc).

\subsection{Regulation by two repressors arranged in an `OR gate'}

Now, we consider regulation of gene $G$ by two independent and repressive TFs, $A$ and $B$. In an OR gate configuration, the presence of one TF is sufficient for the repression of $G$. Similar to the previous section, we have
\be
A + e_{A} \;\autorightleftharpoons{$\alpha_0$}{$\alpha_1$} \;r_{A}
\qquad {\rm and} \qquad
B + e_{B} \;\autorightleftharpoons{$\beta_0$}{$\beta_1$} \;r_{B}
\label{procOR}
\ee
with steady states
\be
[e_A] = \frac{\alpha_1}{\alpha_1 + \alpha_0 [A]} = \frac{1}{1 + \alpha [A]}
\qquad {\rm and} \qquad
[e_B] = \frac{\beta_1}{\beta_1 + \beta_0 [B]} = \frac{1}{1 + \beta [B]}. \label{rarb}
\ee
where $\alpha = \alpha_0/\alpha_1$ and $\beta = \beta_0/\beta_1$.

Because $A$ and $B$ are repressing and act independently, gene $G$ remains active when both binding sites are empty. Accordingly,
\be
[G^*] = [e_A] [e_B] = (1 - r_A) (1 - r_B).
\ee
Therefore, the dynamics of the protein $P$ is described by
\be
\frac{d[P]}{dt} = g \frac{1}{1 + \alpha [A]} \cdot
\frac{1}{1 + \beta [B]} - d_P [P].
\ee

\subsection{Regulation by two repressors arranged in an `AND gate'}

In an AND gate configuration, the presence of both TFs $A$ and $B$ is necessary for the repression of $G$. Accordingly, gene $G$ is active when at least one of the binding sites is empty. In this case, the average number of activated genes is
\be
[G^*] = [G_1^*] + [G_2^*] + [G_3^*] = [e_A] [e_B] + [r_A] [e_B] + [e_A] [r_B].
\ee
where $[r_A]$ and $[r_B]$ are given by equation (\ref{rarb}). For generality, we here assume that the production rate of protein $P$ depends on the configuration of bound TF \cite{Andrecut11}. Therefore, by assuming  production rates $g_1$, $g_2$ and $g_3$ for $[G_1^*]$, $[G_2^*]$ and $[G_3^*]$, respectively, the dynamics of protein $P$ is described by
\be
\frac{d[P]}{dt} = g_1 \frac{1}{1 + \alpha[A]} \cdot \frac{1}{1 + \beta[B]} + g_2 \frac{\alpha[A]}{1 + \alpha[A]} \cdot \frac{1}{1 + \beta[B]} + g_3 \frac{1}{1 + \alpha[A]} \cdot \frac{\beta[B]}{1 + \beta[B]} - d_P [P].
\ee
This can be rewritten as
\be
\frac{d[P]}{dt} = \frac{g_1 + g_2 \alpha [A] + g_3 \beta [B]}{(1 + \alpha[A])(1 + \beta[B])} - d_P [P],
\ee
which is the form used in (3) and (5).

\subsection{Regulation by three repressors arranged in an `AND gate'}

Similarly to the AND gate with two repressors, derived above, the average number of activated genes for the AND gate with three repressors is 
\be
[G^*] = [e_A] [e_B] [e_C] + [e_A] [e_B] [r_C] + [e_A] [r_B] [e_C] + [r_A] [e_B] [e_C] + [e_A] [r_B] [r_C] + [r_A] [e_B] [r_C] + [r_A] [r_B] [e_C].
\ee
Accordingly, assuming different maximal expression rate for each configuration of binding gives
\be
\frac{d[P]}{dt} = \frac{g_1 + g_2 \alpha [A] + g_3 \beta [B] + g_4 \gamma [C] + g_5 \beta\gamma [B][C] + g_6 \alpha\gamma [A][C] + g_7 \alpha\beta [A][B]}{(1 + \alpha[A])(1 + \beta[B])(1 + \gamma[C])} - d_P [P],
\ee
which is the form used in (4).

\section{Analysis of the oscillation period in the standard repressilator}

\begin{figure}[ht]
\begin{center}
\includegraphics[height=6.0cm, angle=0]{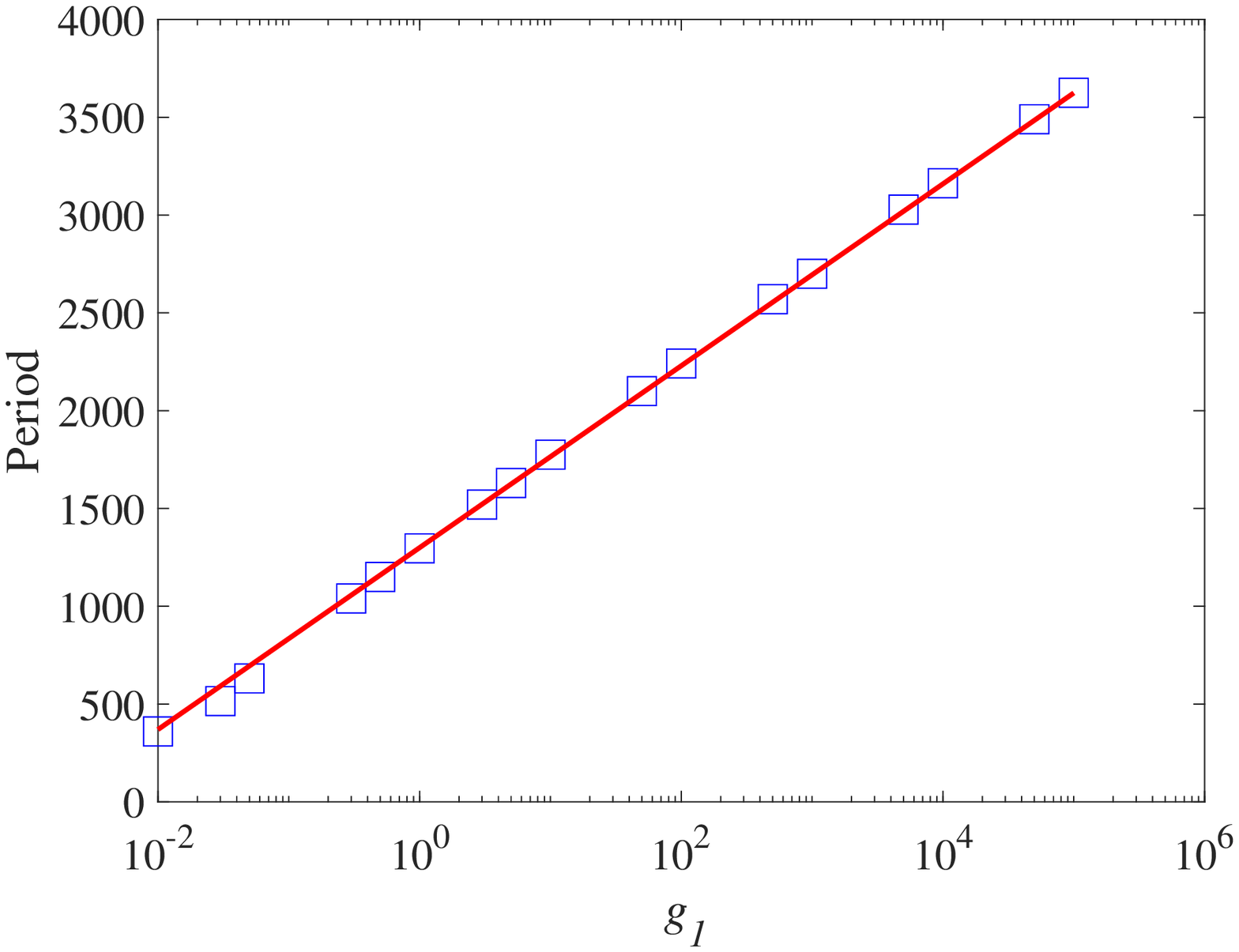}
\includegraphics[height=6.0cm, angle=0]{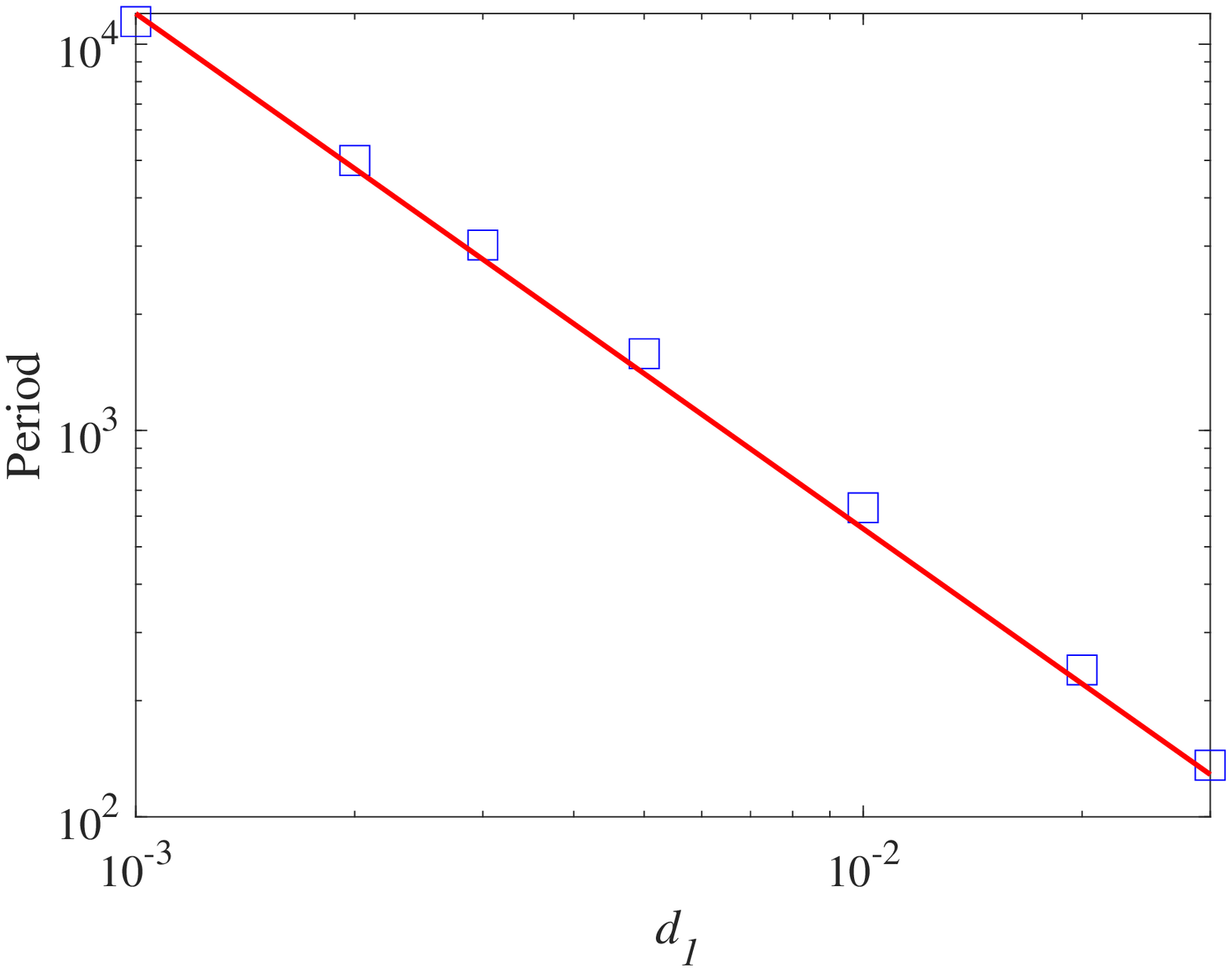}
\end{center}
\caption{{\bf Fit of the oscillations period dependency on $g_1$ and $d_1$ parameters in the repressilator model.} In both panels the points (blue squares) represent the result of the simulations for $d_1=d_2=d_3=0.01$ fixed (left) and $g_1=g_2=g_3=0.05$ fixed (right), while the solid lines correspond to best fitting curves as detailed in the text. The resulting fitting parameters are $a_0=1300$, $a_1=202$, $a_2=1.2$, and $\beta=4/3$. In all cases $b_1=b_2=b_3=10^{-5}$ and $\alpha_1=\alpha_2=\alpha_3=50$.}
\label{fig:fig16}
\end{figure}

Here we show a numerical exploration of the dependency of the period of oscillations of the standard repressilator on the two significant parameters in the model. We restrict attention to the fully-symmetric version for simplicity: $g_1=g_2=g_3$, $d_1=d_2=d_3$, $\alpha_1=\alpha_2=\alpha_3$ and $b_1=b_2=b_3$ in all cases. We consider separately the variation in the period $T$ of the oscillation with
changes in $g_1$ and in $d_1$; these parameters are as defined in (1).
Our simulations are presented in \cref{fig:fig16}. 

The fitting functions assumed in these cases are as follows:
\be
T(g_1) = a_0 + a_1 \log(g_1) \qquad {\rm for} \;d_1 \; {\rm fixed} 
\ee
and 
\be
T(d_1) = a_2 d_1^{-\beta} \qquad {\rm for} \;g_1 \; {\rm fixed} 
\ee
The simulations show a very weak dependency of the period on the parameter $g_1$ - over seven
orders of magnitude in $g_1$ the period increases by a factor of only around seven.
As $d_1$ increases the period drops much more rapidly, in line with the power-law
scaling $T \propto d_1^{-4/3}$. This is to be expected since the degradation rate
$d_1$ is a key
timescale in the problem and the decay rate of the concentrations of TFs is essential to the
nature of the `repressilator' dynamics.
  In neither case does the qualitative nature of the periodic
oscillation change over these ranges of parameters. This allows us to conclude that the
parameter values used, for example in figure 2 in the main text, are entirely typical of 
the repressilator's dynamical behaviour. Finally we note that the dependence of the period
and shape of the oscillation on the production
rates $b_1,b_2$ and $b_3$ is very weak in this parameter range where typical concentrations of
TFs are maintained well above each of the $b_i$ by the nonlinear prodution terms.

\section{Analysis of mathematical concepts in the cross-repressilator}

Three mathematical concepts that need to be explained and analysed more closely are symmetric structure of systems~(2) and (3), the growth rate of period of limit cycles close to SNIC/heteroclinic bifurcation and determining factors in the order of expressed genes when the system oscillates.

\subsection{Analysis of the three-gene Models in the symmetric case}

In the case that the models~(2) and (3) are cyclically symmetric under the permutation $(x_1,x_2,x_3) \rightarrow (x_2,x_3,x_1)$, i.e. when the parameters $g$, $\alpha$, $\beta$, $k$ and $h$ do not take different values in each equation, the symmetric structure of the problem enables some semi-analytic progress to be made since the Jacobian matrix at the symmetric equilibrium point $(x,x,x)$ is constrained to always take the form
\be
Df|_{(x,x,x)} = \left(
\begin{array}{ccc}
a & b & c \\
c & a & b \\
b & c & a
\end{array}
\right), \label{eqn:jacobian}
\ee
where the matrix entries are the partial derivatives of the function on the right-hand side, say $\dot x_1 = f_1(x_1,x_2,x_3)$. In particular we have
\be
a \equiv \left. \frac{\partial f_1}{\partial x_1}\right|_{(x,x,x)}, \qquad
b \equiv \left. \frac{\partial f_1}{\partial x_2}\right|_{(x,x,x)}, \quad \mathrm{and} \quad
c \equiv \left. \frac{\partial f_1}{\partial x_3}\right|_{(x,x,x)},
\ee
where the partial derivatives are evaluated at the symmetric equilibrium point which can usually not be solved for explicitly.
The eigenvalues of~\cref{eqn:jacobian} are then a single real eigenvalue and a complex conjugate pair
\begin{align}
\lambda_1 &= a+b+c, \\
\lambda_{2,3} &= a - \frac{1}{2}(b+c) \pm \frac{i\sqrt{3}}{2}(b-c).
\end{align}

In the case of the original repressilator model~(1) for which $f_1(x_1,x_2,x_3)=b+\frac{g}{1+\alpha x_2^h}-k x_1$, we see that the location of the equilibrium point satisfies $\frac{g}{1+\alpha x^h}=kx-b$. Its eigenvalues are therefore
\begin{align}
\lambda_1 &= -k -\frac{\alpha h}{g}x^{h-1} (kx-b)^2, \\
\lambda_{2,3} &= -k + \frac{\alpha h}{2g}x^{h-1} (kx-b)^2 \pm \frac{i\sqrt{3}\alpha h}{2g}x^{h-1}(kx-b)^2.
\end{align}
So $\lambda_1<0$ always, and at the Hopf bifurcation point, where $Re(\lambda_{2,3}) =0$ we see that the frequency
$Im(\lambda_{2,3}) = k \sqrt{3}$ independent of the other parameters. In the symmetric case this implies that the frequency
at the two Hopf bifurcations from the symmetric equilibrium, illustrated in~Fig.~4~(left), would be equal.

For the cross-repressilator with the `OR' gate we have $f_1(x_1,x_2,x_3)=\frac{g}{(1+\alpha x_2^h)(1+\beta x_3^h)} - k x_1$. Hence the symmetric equilibrium $(x,x,x)$ satisfies $(1+\alpha x^h)(1+\beta x^h)=\frac{g}{kx}$. Computing the relevant partial derivatives yields the three eigenvalues
\begin{align}
\lambda_1 &= -k -\frac{k\alpha h x^h}{1+\alpha x^h} - \frac{k\beta h x^h}{1+\beta x^h} \\
\lambda_{2,3} &= -k +\frac{1}{2}\left( \frac{k\alpha h x^h}{1+\alpha x^h} + \frac{k\beta h x^h}{1+\beta x^h} \right)
\pm \frac{i\sqrt{3}k^2 h(\beta-\alpha)x^{1+h}}{2g}
\end{align}
which shows that in this case also, $\lambda_1<0$ and there is a possibility of a Hopf bifurcation where $Re(\lambda_{2,3})=0$.
The imaginary part of the complex conjugate pair is proportional to $(\alpha-\beta)$ which indicates, in agreement with
numerical explorations, that these eigenvalues become real along the line $\alpha=\beta$.

Finally we note that the occurrence of the SNIC bifurcation and the codimension-two point indicated by the grey dot in~Fig.~6 are also consequences of symmetry. We do not comment further here, and leave a more detailed analysis of the consequences of symmetry for this problem to be the subject of a future paper \cite{Dawes2021}.

\subsection{Analysis of the oscillation period close to a global (SNIC/heteroclinic) bifurcation}

\begin{figure}[ht]
\begin{center}
\includegraphics{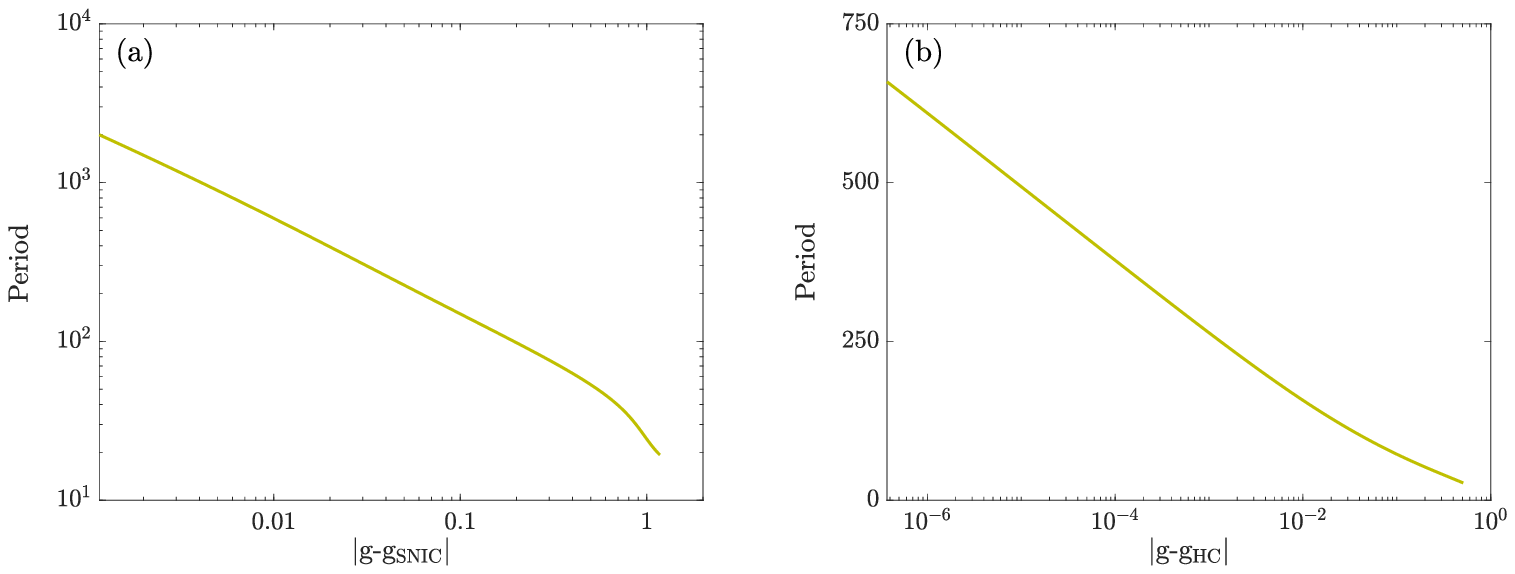}
\end{center}
\caption{{\bf Growth scale of periodic orbits.} Period of the family of periodic orbits ending at a SNIC bifurcation for $\alpha=9$ (left) and heteroclinic bifurcation for $\alpha=1$ (right). The rest of parameters of the system are the same as~Fig.~4(b).}
\label{fig:fig14}
\end{figure}

In this subsection we investigate in more detail the scaling behaviour of the increase of the period of the limit cycles as they terminate at two different global bifurcations: the saddle-node on an invariant circle (SNIC) bifurcation and the heteroclinic bifurcation. Appropriate plots of the period of the orbit near these global bifurcations yields straight-line plots. For both of these bifurcations, the period increases in accordance with a well-known analytical scaling law when $g\longrightarrow g_{\rm SNIC}$. Close to a SNIC bifurcation, it is known that $T\propto 1/\sqrt{|g-g_{\rm SNIC}|}$, i.e. a straight line in log-log scale. The left panel in \cref{fig:fig14} shows a transformation of~Fig.~4(b) in log-log scale which therefore matches with what we expect for a SNIC bifurcation. For $\alpha=9$ we find $g_{\rm SNIC}\approx 1.51449$.

When the periodic orbit terminates at a heteroclinic bifurcation, the period increases according to $T\propto \log{|g-g_{\rm HC}|}$ as $g\longrightarrow g_{\rm HC}$. Therefore a plot of the period $T$ against $g-g_{\rm HC}$ on a log-linear scale should be a straight line. We confirm this in the right panel of \cref{fig:fig14} where we show the data computed for~Fig.~6~(right) but on a log-linear scale, again shifted by subtraction of $g_{\rm HC}\approx 1.05509$.

A further detailed verification of our numerical work is that the theoretical result for the constant in the relation $T\propto \log{1/|g-g_{\rm HC}|}$ is given by 
$3/\lambda_+$, i.e. $T \sim (3/\lambda_+) \log(1/|g-g_{\rm HC}|)$ where $\lambda_+$ is the positive eigenvalue of the saddle point equilibrium that collides with the periodic orbit (in the case that the periodic orbit is stable, which it is here). The factor of 3 arises due to the $\mathbb{Z}_3$ symmetry and the three distinct saddle points that all collide with the periodic orbit simultaneously. Numerically, we find that $\lambda_+ \approx 0.05922$ (5 dp); therefore $1/\lambda_+ = 16.88619$ (5 dp).  Hence the theoretical estimate for the slope of the straight line is $3\times 16.88619 = 50.65857$ (5 dp). From the right panel of \cref{fig:fig14} we estimate the slope of the curve to be $50.17112$ based on the values of the period T at $|g-g_{\rm HC}|= 3\times 10^{-7}$ and $|g-g_{\rm HC}|\approx 3\times 10^{-4}$. This therefore provides quantitative confirmation of the bifurcation behaviour in the differential equations.

\subsection{Ordering of the gene expression in the oscillations}

As discussed in the main text, the genes are expressed in a particular order in the oscillatory regimes. For instance, in the three-gene case with the OR gate configuration, the gene expression occurs cyclically clockwise. That is, specifically, first gene $\sf X_3$ is expressed and then genes $\sf X_2$ and finally $\sf X_1$. This order of gene expression is due to the values of the parameters of~(2). In~(2), we chose $\alpha=9>\beta=0.1$ which makes the inhibitions of the inner loop in~Fig.~3 stronger; this results in the genes being expressed in clockwise order. However, if we set $\alpha<\beta$, the order of gene expression in the oscillatory regime is reversed.

\begin{figure}[ht]
\begin{center}
\includegraphics{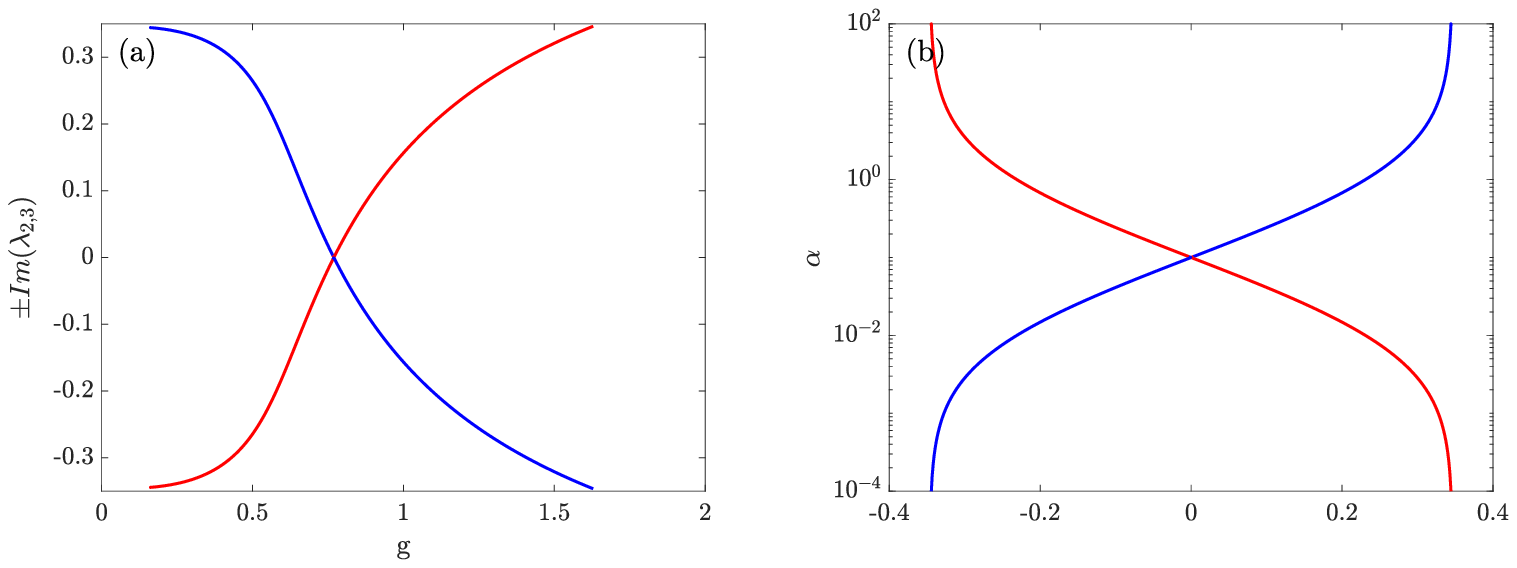}
\caption{{\bf Imaginary parts of eigenvalues along the Hopf curve.} The variation in the imaginary parts of the complex conjugate eigenvalues along the Hopf bifurcation curve ($\sf HB$, dark green) in~Fig.~6 shown versus $g$ (left panel) and $\alpha$ in log scale (right panel). Other parameter values are as in~Fig.~6.}
\label{fig:fig15}
\end{center}
\end{figure}

Mathematically the oscillatory direction of the family of periodic orbits emanating from a (supercritical) Hopf bifurcation is determined by the relative signs of the real and imaginary parts of the eigenvectors corresponding to the pair of complex conjugate eigenvalues ($\lambda_{2,3}=Re(\lambda_{2,3})\pm Im(\lambda_{2,3})$). Since these eigenvectors remain non-zero and vary continuously along the Hopf bifurcation curve on which these eigenvalues have zero real part, we can conclude that the oscillation changes direction when the imaginary parts of the eigenvalues pass through zero. \Cref{fig:fig15} shows the imaginary parts of $\lambda_{2,3}$ and its conjugate $-\lambda_{2,3}$ along the curve of Hopf bifurcations (dark green) in~Fig.~6 projected on $g$ and $\alpha$ (in log scale) axes. For small values of $g$ (large values of $\alpha$), let assume that $Im(\lambda_{2,3})$ is positive (blue); therefore, its conjugate, $-Im(\lambda_{2,3})$, is negative (red). As $g$ increases ($\alpha$ decreases), $Im(\lambda_{2,3})$ and its conjugate, $-Im(\lambda_{2,3})$, cross each other. Thus, the direction of oscillation of periodic orbits switches. The crossing point of $Im(\lambda_{2,3})$ and $-Im(\lambda_{2,3})$ at 0 corresponds to the point that the Hopf curve ($\sf HB$) touches the homoclinic curve ($\sf HC$, red). Precisely at the parameter value where the complex conjugate pair of eigenvalues are zero, the Hopf bifurcation theorem breaks down and further analysis is required; this discussion will be presented in future work \cite{Dawes2021}.

\section{Numerical Methods}

We have employed different tools and software packages for computing solutions, bifurcation diagrams and colourmaps. All the time series have been integrated in XPPAUT, freely available at \url{http://www.math.pitt.edu/~bard/xpp/xpp.html}. The bifurcation diagrams are computed using pseudo-arclength continuation methods in software package {\sc Auto} \cite{Doedel81,Doedel10}.

Finally, the colourmaps are computed in MATLAB$\textsuperscript{\textregistered}$ (Mathworks Inc.) starting from different initial conditions. In~Fig.~5, the initial-condition planes are colour-coded depending on the steady state (one of the three) that the system converges to. Similarly, in~Fig.~8, the final state of the system when $g$ increases and saturates at a value larger than $g_{\rm SNIC}$ determines the colour of the initial conditions.